\documentclass[11pt,leqno]{amsart}
\usepackage[T1]{fontenc}
\usepackage{ae}
\usepackage{aecompl}
\usepackage[cm]{aeguill}
\usepackage{amsmath}
\usepackage{amssymb}
\usepackage{latexsym,graphicx,chicago}
\usepackage{graphicx,multirow,rotating}

\newtheorem{thm}{Theorem}[section]
\newtheorem{lem}{Lemma}[section]

\newtheorem{cor}{Corollary}[section]

\newtheorem{rem}{Remark}[section]

\def\1{{{\mbox{${\rm{1\negthinspace\negthinspace I}}$}}}}

\newcounter{hypc}

\newcounter{cond}

\begin{document}
\title[CLT for unbounded functions of the intermittent map]{Some unbounded functions of  intermittent maps for which the central limit theorem holds}
\author{J. Dedecker$^1$}\thanks{$^1$ Universit\'e Paris 6, Laboratoire de
Statistique Th\'eorique et Appliqu\'ee.}
\address{J. Dedecker\\
 Laboratoire de Statistique Th\'eorique et Appliqu\'ee \newline
Universit\'e paris 6, 175 rue du Chevaleret  75013 Paris, France.
\newline email:
dedecker@ccr.jussieu.fr} \author{C. Prieur$^2$}\thanks{$^2$ INSA
Toulouse, Institut Math\'ematique de Toulouse.} \address{C. Prieur \\
INSA Toulouse, Institut Math\'ematique de Toulouse \newline \'Equipe
de Statistique et Probabilit\'es, 135 avenue de Rangueil, 31077
Toulouse Cedex 4, France.
\newline email: Clementine.Prieur@insa-toulouse.fr}  \maketitle

\noindent{\bf Abstract.} We compute some dependence coefficients for
the stationary Markov chain whose transition kernel is the
Perron-Frobenius operator of an expanding map $T$  of $[0, 1]$ with
a neutral fixed point. We use these coefficients to prove a central
limit theorem for the partial sums of  $f\circ T^i$, when $f$
belongs to a large class of unbounded functions from $[0, 1]$ to
${\mathbb R}$. We also prove  other limit theorems and moment
inequalities.

\medskip

 \noindent {\bf Classification MSC 2000.} 37E05, 37C30, 60F05.

\medskip



\bibliographystyle{chicago}
\section{Introduction}
For  $\gamma$ in $]0, 1[$, we consider the intermittent map
$T_\gamma$ from $[0, 1]$ to $[0, 1]$, studied for instance  by
Liverani, Saussol and Vaienti (1999), which is a modification of the
Pomeau-Manneville map (1980):
$$
   T_\gamma(x)=
  \begin{cases}
  x(1+ 2^\gamma x^\gamma) \quad  \text{ if $x \in [0, 1/2[$}\\
  2x-1 \quad \quad \quad \ \  \text{if $x \in [1/2, 1]$}
  \end{cases}
$$
We denote by $\nu_\gamma$ the unique $T_\gamma$-probability measure
on $[0, 1]$. We denote by $K_\gamma$ the Perron-Frobenius operator
of $T_\gamma$ with respect to $\nu_\gamma$: for any bounded
measurable functions $f, g$,
$$
\nu_\gamma(f \cdot  g\circ T_\gamma)=\nu_\gamma(K_\gamma(f) g) \, .
$$
Let $(X_i)_{i \geq 0}$ be a stationary Markov chain with invariant
measure $\nu_\gamma$ and transition Kernel $K_\gamma$. It is well
known (see for instance Lemma XI.3 in Hennion and Hervé (2001)) that
on the probability space $([0, 1], \nu_\gamma)$, the random variable
$(T_\gamma, T^2_\gamma, \ldots , T^n_\gamma)$ is distributed as
$(X_n,X_{n-1}, \ldots, X_1)$. Hence any information on the law of
$$
S_n(f)=\sum_{i=1}^n f \circ T_\gamma^i
$$
can be obtained by studying  the law of  $\sum_{i=1}^n f(X_i)$.

In 1999, Young proved that such systems (among many others) may be
described  by a Young tower with polynomial decay of the return
time. From this construction, she was able to control the
covariances $\nu_\gamma(f\circ T^n \cdot (g-\nu_\gamma(g)))$ for any
bounded function $f$ and any $\alpha$-H\"older function $g$, and
then to prove that $n^{-1/2}(S_n(f)-\nu_\gamma(f))$ converges in
distribution to a normal law as soon as  $\gamma<1/2$ and $f$ is any
$\alpha$-H\"older function. For $\gamma=1/2$, Gou\"ezel (2004)
proved that  the central limit theorem remains true with the same
normalization $\sqrt n$  if $f(0)=\nu_\gamma(f)$, and with the
normalization $\sqrt{n \ln(n)}$ if $f(0)\neq \nu_\gamma(f)$. When
$1/2<\gamma<1$, he proved that if $f$ is $\alpha$-H\"older and
$f(0)\neq \nu_\gamma(f)$, $n^{-\gamma}(S_n(f)-\nu_\gamma(f))$
converges to a stable law.

At this point, two questions (at least) arise: 1) what happens if
$f$ is no longer continuous? 2) what happens if $f$ is no longer
bounded? For instance, for the uniformly expanding map
$T_0(x)=2x-[2x]$, the central limit theorem holds with the
normalization $\sqrt n$ as soon as $f$ is monotonic and square
integrable on $[0, 1]$, that is not necessarily continuous nor
bounded.

For the slightly different map
$\theta_\gamma(x)=x(1-x^\gamma)^{-1/\gamma}-[x(1-x^\gamma)^{-1/\gamma}]$,
with the same behavior around the indifferent fixed point, Raugi
(2004) (following a work by Conze and Raugi (2003)) has given a
precise criterion for the central limit theorem with the
normalization $\sqrt n$ in the case where $0<\gamma<1/2$ (see his
Corollary 1.7). In particular his result applies to a large class of
non continuous functions, which gives a quite complete answer to our
first question for the map $\theta_\gamma$. The result also applies
to the unbounded function $f(x)=x^{-a}$ with $0<a<1/2-\gamma$.
However, the function $f$ is allowed to blow up near $0$ only (if
$f$ tends to infinity when $x$ tends to $x_0\in ]0, 1]$, then the
variation coefficient $v(fh_\gamma,k)$, where $h_\gamma$ is the
density of the $\theta_\gamma$-invariant probability, is always
infinite).

We now go back to the map $T_\gamma$. In a short discussion after
the proof of his Theorem 1.3, Gou\"ezel (2004) considers the case
where $f(x)=x^{-a}$, with $0<a<1-\gamma$. He shows that, if
$0<a<1/2-\gamma$ then the central limit theorem holds with the
normalization $\sqrt n$, if $a=1/2-\gamma$ then the central limit
theorem holds with the normalization $\sqrt{n \ln(n)}$, and if
$0<a<1-\gamma$ and $\gamma\geq 1/2$ then there is convergence to a
stable law. Again, as for Raugi's result (2004) concerning the map
$\theta_\gamma$, the function $f$ is allowed to blow up only near
$0$.


On another hand, we know that for stationary Harris recurrent Markov
chains with invariant measure $\mu$ and $\beta$-mixing coefficients
of order $n^{-b}$, $b>1$,   the central limit theorem holds with the
normalization $\sqrt n$ as soon as the  moment condition
$\mu(|f|^p)<\infty$ holds for $p>2b/(b-1)$. For $T_\gamma$, the
covariances decay is of order $n^{(\gamma-1)/\gamma}$, so that one
can expect the moment condition $\nu_\gamma(|f|^p)<\infty$  for
$p>(2-2\gamma)/(1-2\gamma)$. For instance, if $f(x)=x^{-a}$, since
the density of $\nu_\gamma$ is of order $x^{-\gamma}$ near $0$,  the
moment condition is satisfied if $0<a<1/2-\gamma$, which is coherent
with Gou\"ezel's result (2004). However, since the chain $(K_\gamma,
\nu_\gamma)$ is not $\beta$-mixing, the condition
$\nu_\gamma(|f|^p)<\infty$ for $p>(2-2\gamma)/(1-2\gamma)$ alone  is
not sufficient to imply the central limit theorem, and one still
needs some regularity on $f$.

Let us now define the class of functions of interest. For any
probability measure $\mu$ on ${\mathbb R}$, any $M>0$ and any $p \in
]1, \infty]$, let $\text{Mon}(M, p, \mu)$ be the class of functions
$g$ which are monotonic on some open interval of ${\mathbb R}$ and
null elsewhere, and such that $\mu(|g|>t)\leq M^p t^{-p}$ for
$p<\infty$ and $\mu(|g|>M)=0$ for $p=\infty$. Let ${\mathcal C}(M,
p, \mu)$ be the closure in ${\mathbb L}^1(\mu)$ of the set of
functions which can be written as $\sum_{i=1}^n a_i g_i$, where
$\sum_{i=1}^n |a_i| \leq 1$ and $g_i$ belongs to $\text{Mon}(M, p,
\mu)$. Note that a function belonging to ${\mathcal C}(M, p, \mu)$
is allowed to blow up at an infinite number of points.

In Corollary \ref{CLTT} of the present paper, we prove that if $f$
belongs to the class ${\mathcal C}(M, p, \nu_\gamma)$ for
$p>(2-2\gamma)/(1-2\gamma)$, then $n^{-1/2}(S_n(f-\nu_\gamma(f))$
converges in distribution to a normal law. We also give some
conditions on $p$ to obtain rates of convergence in the central
limit theorem (Corollary \ref{Krov}), as well as moment inequalities
for $S_n(f-\nu_\gamma(f))$ (Corollary \ref{melnic}). Finally, a
central limit theorem for the empirical distribution function of
$(T_\gamma^i)_{1\leq i \leq n}$ is given in the last section
(Corollary \ref{empT}).

To prove these results, we compute the $\beta$-dependence
coefficients (cf Dedecker and Prieur (2005, 2007)) of  the Markov
chain $(K_\gamma, \nu_\gamma)$. The main tool is a precise estimate
of the Perron-Frobenius operator  of the map $F$ associated to
$T_\gamma$ on the Young tower, due to Maume-Deschamps (2001). Next,
we apply some general results for $\beta$-dependent Markov chains.
For the sake of  simplicity, we give all the computations in the
case of the maps $T_\gamma$, but our arguments remain valid for many
other systems modelled by Young towers.

\section{The main inequality}
For any Markov kernel $K$ with invariant measure $\mu$, any
non-negative integers $n_1, n_2, \ldots , n_k$, and any bounded
measurable functions $f_1, f_2, \ldots, f_k$, define
\begin{eqnarray*}
K^{(n_1, n_2, \ldots, n_k)}(f_1, f_2, \ldots, f_k)&=&K^{n_1}(f_1
K^{n_2}(f_2 K^{n_3}(f_3\cdots K^{n_{k-1}}(f_{k-1}K^{n_k}(f_k))\cdots
)))\, , \ \text{and}\\
K^{(0)(n_1, n_2, \ldots, n_k)}(f_1, f_2, \ldots, f_k)&=&K^{(n_1,
n_2, \ldots, n_k)}(f_1, f_2, \ldots, f_k)-\mu(K^{(n_1, n_2, \ldots,
n_k)}(f_1, f_2, \ldots, f_k))\, .
\end{eqnarray*}
For $\alpha \in ]0, 1]$ and $c>0$, let $H_{\alpha, c}$ be the set of
functions $f$ such that $|f(x)-f(y)|\leq c|x-y|^\alpha$.
\begin{thm}\label{maum}
Let $\gamma \in ]0, 1[$, and let $f^{(0)}=f-\nu_\gamma (f)$. For any
$\alpha \in ]0, 1]$, the following inequality holds:
$$
\nu_\gamma\Big(\sup_{f_1, \ldots , f_k \in H_{\alpha, 1}}\big|
K_\gamma^{(0)(n_1, n_2, \ldots, n_k)}(f_1^{(0)}, f_2^{(0)}, \ldots,
f_k^{(0)})\big| \Big) \leq \frac{C(\alpha,
k)(\ln(n_1+1))^{2}}{(n_1+1)^{(1-\gamma)/\gamma}}\, .
$$
In particular,
$$
\nu_\gamma\Big(\sup_{f \in H_{\alpha
,1}}|K_\gamma^nf-\nu_\gamma(f)|\Big)\leq
\frac{C(\alpha,1)(\ln(n+1))^{2}}{(n+1)^{(1-\gamma)/\gamma}}\, .
$$
\end{thm}
\noindent {\bf Proof of Theorem \ref{maum}.}  We refer to the paper
by Young (1999) for the construction of the tower $\Delta$
associated to $T_\gamma$ (with floors $\Lambda_\ell$), and for the
mappings $\pi$ from $\Delta$ to $[0, 1]$ and $F$ from $\Delta$ to
$\Delta$ such that $T_\gamma\circ \pi=\pi \circ F$. On $\Delta$
there is a probability measure $m_0$ and an unique $F$-invariant
probability measure $\bar \nu$ with density $h_0$ with respect to
$m_0$, and $\bar \nu (\Lambda_\ell)=O(\ell^{-1/\gamma})$. The unique
$T_\gamma$-invariant probability measure $\nu_\gamma$ is then given
by $\nu_\gamma=\bar \nu^\pi$. There exists a distance $\delta$ on
$\Delta$ such that $\delta(x, y)\leq 1$ and $|\pi(x)-\pi(y)|\leq
\kappa \delta(x,y)$. For $\alpha \in ]0, 1]$,  let
$\delta_\alpha=\delta^\alpha$, let $L_\alpha$ be the space of
Lipschitz functions with respect to $\delta_\alpha$, and let
$L_\alpha(f)=\sup_{x,y \in \Delta}|f(x)-f(y)|/\delta_\alpha(x,y)$.
Let $L_{\alpha,c}$ be the set of functions such that
$L_\alpha(f)\leq c$. For $\varphi$ in $H_{\alpha, c}$, the function
$\varphi\circ \pi$ belongs to $L_{\alpha,c\kappa^\alpha}$. Any
function $f$ in $L_\alpha$ is bounded and the space $L_\alpha$ is a
Banach space with respect to the norm
$\|f\|_\alpha=L_\alpha(f)+\|f\|_\infty$. The density $h_0$ belongs
to any $L_\alpha$ and $1/h_0$ is bounded. As in Maume-Deschamps
(2001), we denote by ${\mathcal L}_0$ the Perron-Frobenius operator
of $F$ with respect to $m_0$, and by $P$ the Perron-Frobenius
operator of $F$ with respect to $\bar \nu$: for any bounded
measurable functions $\varphi, \psi$,
$$
 m_0  (\varphi \cdot \psi \circ F)=m_0({\mathcal L}_0(\varphi) \psi) \quad
 \text{and} \quad
\bar \nu  (\varphi \cdot \psi \circ F)=\bar \nu (P(\varphi) \psi)\,
.
$$

We first state a useful lemma
\begin{lem}\label{deb} For any positive  $n_1, n_2, \ldots , n_k$
and any bounded measurable functions $f_{1}, f_{2}, \ldots ,f_{k}$
from $[0, 1]$ to ${\mathbb R}$, one has $$\label{cqf}
K_\gamma^{(n_1, n_2, \ldots, n_k)}(f_1, f_2, \ldots, f_k) \circ
\pi={\mathbb E}_{\bar \nu}\big(P^{(n_1, n_2, \ldots, n_k)}(f_1\circ
\pi, f_2\circ \pi, \ldots, f_k\circ \pi) \big|\pi\big)\, .
$$
\end{lem}
We now complete the proof of Theorem \ref{maum} for $k=2$, the
general case being similar. Applying Lemma \ref{deb}, it follows
that \begin{multline*} \sup_{f, g \in H_{\alpha,
1}}|K_\gamma^n(f^{(0)}K_\gamma^m
g^{(0)})(x)-\nu_\gamma(f^{(0)}K_\gamma^m g^{(0)})|\\ \leq {\mathbb
E}_{\bar \nu}\Big(\sup_{\phi, \psi \in L_{\alpha,\kappa^\alpha}}
|P^n(\phi^{(0)}P^m\psi^{(0)})-\bar \nu (\phi^{(0)}P^m
\psi^{(0)})|\Big |\pi=x\Big)\, .
\end{multline*}
Here, we need the following lemma, which is derived from Lemma 3.4
in Maume-Deschamps (2001).
\begin{lem}\label{maum2} There exists $M_\alpha>0$
such that, for any $\psi \in L_\alpha$, $$|P^m
\psi(x)-P^m\psi(y)|\leq M_\alpha
\delta(x,y)\|\psi^{(0)}\|_\alpha\leq 2M_\alpha
\delta_\alpha(x,y)L_\alpha(\psi)\, .$$
\end{lem} Hence, if $\psi \in  L_{\alpha, \kappa^\alpha}$, then  $P^m(\psi^{(0)})$ belongs to $L_{\alpha, 2M_\alpha \kappa^\alpha}$ and is
centered, so that $\phi^{(0)}P^m \psi^{(0)}$ belongs to $L_{\alpha,
4M_\alpha\kappa^{2\alpha}}$. It follows that
$$
\sup_{f, g \in H_{\alpha, 1}}|K_\gamma^n(f^{(0)}K_\gamma^m
g^{(0)})(x)-\nu(f^{(0)}K_\gamma^m g^{(0)})|\leq
4M_\alpha\kappa^{2\alpha}{\mathbb E}_{\bar \nu}\Big(\sup_{\varphi
\in L_{\alpha, 1}} |P^n(\varphi)-\bar \nu (\varphi)|\Big
|\pi=x\Big)\, .
$$
Next, we apply the following Lemma, which is derived from Corollary
3.14 in Maume-Deschamps (2001).
\begin{lem}\label{maum3} Let
$v_\ell=(\ell+1)^{(1-\gamma)/\gamma}(\ln(\ell+1))^{-2}$. There
exists $C_\alpha>0$ such that
$$
{\mathbb E}_{\bar \nu}\Big(\sup_{\varphi \in L_{\alpha,1}}
|P^n(\varphi)-\bar \nu (\varphi)|\Big |\pi=x\Big)\leq
C_\alpha(\ln(n+1))^{2}(n+1)^{(\gamma-1)/\gamma}\sum_{\ell \geq
0}v_\ell {\mathbb E}_{\bar \nu}({\bf 1}_{\Lambda_\ell}|\pi=x)\, .
$$
\end{lem}
Hence
$$
\nu_\gamma\Big(\sup_{f, g \in H_{\alpha,
1}}|K_\gamma^n(f^{(0)}K_\gamma^m g^{(0)})-\nu(f^{(0)}K_\gamma^m
g^{(0)})|\Big)\leq 4M_\alpha \kappa^{2\alpha}
C_\alpha(\ln(n+1))^{2}(n+1)^{(\gamma-1)/\gamma}\sum_{\ell \geq
0}v_\ell \bar \nu(\Lambda_\ell)\, .
$$
Since $\bar \nu(\Lambda_\ell)=O(\ell^{-1/\gamma})$, the result
follows.

\medskip

\noindent {\bf Proof of Lemma \ref{deb}.} We write the proof for
$k=2$ only, the general case being similar. Let $\varphi, f$ and $g$
be three bounded measurable functions. One has
\begin{eqnarray*}
\nu_\gamma( \varphi K_\gamma^n (fK_\gamma^mg))&=&\nu_\gamma( \varphi \circ T_\gamma^{n+m} \cdot f\circ T_\gamma^m \cdot g)\\
&=& \bar \nu (\varphi \circ \pi \circ F^{n+m} \cdot f\circ \pi \circ F^m \cdot g\circ \pi)\\
&=& \bar \nu (\varphi \circ \pi  P^n (f\circ \pi  P^m (g\circ \pi)))\\
 &=& \bar \nu (\varphi \circ \pi {\mathbb E}_{\bar \nu} (P^n(f \circ \pi P^m (g\circ \pi))|\pi))\\
&=& \int \varphi(x){\mathbb E}_{\bar \nu}(P^n(f \circ \pi P^m
(g\circ \pi))|\pi=x)\nu_\gamma(dx)\, ,
\end{eqnarray*}
which proves Lemma \ref{deb} for $k=2$.

\medskip

\noindent {\bf Proof of Lemma \ref{maum2}.} Applying Lemma 3.4 in
Maume-Deschamps (2001) with $v_k=1$, we see that there exists
$D_\alpha>0$ such that, for any $\psi$ in $L_\alpha$,
$$
|{\mathcal L}_0^m \psi (x)-{\mathcal L}_0^m \psi(y)|\leq D_\alpha
\delta_\alpha(x, y)\|\psi\|_\alpha  .
$$
Now $P^m(\psi)={\mathcal L}_0^m(\psi h_0)/h_0$. Since $1/h_0$ is
bounded by $B(h_0)$, and since $h_0$ belongs to $L_\alpha$, it
follows that
$$
|P^m \psi (x)-P^m \psi(y)|\leq D_\alpha B(h_0)\|h_0\|_\alpha
\delta_\alpha (x, y)\|\psi\|_\alpha .
$$
Let $M_\alpha=D_\alpha B(h_0)\|h_0\|_\alpha$. Since $|P^m \psi
(x)-P^m \psi(y)|=|P^m \psi^{(0)} (x)-P^m \psi^{(0)}(y)|$ and since
$\|\psi^{(0)}\|_\infty \leq L_\alpha(\psi)$, it follows that
$$
|P^m \psi (x)-P^m \psi(y)|\leq M_\alpha
\delta_\alpha(x,y)\|\psi^{(0)}\|_\alpha\leq 2M_\alpha \delta_\alpha
(x,y) L_\alpha(\psi)\, .
$$

\medskip

\noindent {\bf Proof of Lemma \ref{maum3}.} Applying Corollary 3.14
in Maume-Deschamps (2001), there exists $B_\alpha>0$ such that
$$
|{\mathcal L}_0^n f-h_0m_0(f)|\leq B_\alpha
\|f\|_\alpha(\ln(n+1))^{2}(n+1)^{(\gamma-1)/\gamma}\sum_{\ell \geq
0} v_\ell {\bf 1}_{\Delta_\ell}\, .
$$
It follows that, with the notations of the proof of Lemma
\ref{maum2},
$$
|P^n (f)-\bar \nu(f)|\leq B_\alpha
B(h_0)\|h_0\|_\alpha\|f\|_\alpha(\ln(n+1))^{2}(n+1)^{(\gamma-1)/\gamma}\sum_{\ell
\geq 0} v_\ell {\bf 1}_{\Delta_\ell}\, .
$$
Since $|P^n (f)-\bar \nu(f)|=|P^n(f^{(0)})-\bar \nu(f^{(0)})|$ and
since $\|f^{(0)}\|_\infty \leq L_\alpha(f)$, it follows that
$$
|P^n (f)-\bar \nu(f)|\leq 2B_\alpha B(h_0)\|h_0\|_\alpha
L_\alpha(f)(\ln(n+1))^{2}(n+1)^{(\gamma-1)/\gamma}\sum_{\ell \geq 0}
v_\ell {\bf 1}_{\Delta_\ell}\, ,
$$
and the result follows.

\section{The dependence coefficients}\label{depcoef}
 Let ${\bf X}=(X_i)_{i\geq 0}$ be a stationary  Markov chain with invariant measure $\mu$ and transition kernel
$K$. Let $f_t(x)={\bf 1}_{x \leq t}$. As in Dedecker and Prieur
(2005, 2007), define the coefficients $\alpha_k(n)$  of the
stationary Markov chain $(X_i)_{i \geq 0}$ by
\begin{eqnarray*}
\alpha_1(n)&=&\sup_{t \in {\mathbb R}}\mu(|K^n(f_t)-\mu(f_t)|)\, , \quad \text{\and for $k\geq 2$,}\\
\alpha_k(n)&=&\alpha_1(n)\vee \sup_{2\leq l \leq k}\sup_{n_2\geq 1,
\ldots n_l\geq 1}\sup_{t_1, \ldots, t_l \in {\mathbb R}}
\mu\big(|K^{(0)(n, n_2, \ldots, n_l)}(f_{t_1}, f_{t_2}, \ldots,
f_{t_l})|\big)\, .
\end{eqnarray*}
In the same way, define the coefficients  $\beta_k(n)$ by
\begin{eqnarray*}
\beta_1(n)&=&\mu\Big(\sup_{t \in {\mathbb R}}|K^n(f_t)-\mu(f_t)|\Big)\, , \quad \text{\and for $k\geq 2$,}\\
\beta_k(n)&=&\beta_1(n)\vee \sup_{2\leq l \leq k}\sup_{n_2\geq 1,
\ldots n_l\geq 1} \mu\Big(\sup_{t_1, \ldots, t_l \in {\mathbb
R}}|K^{(0)(n, n_2, \ldots, n_l)}(f_{t_1}, f_{t_2}, \ldots,
f_{t_l})|\Big)\, .
\end{eqnarray*}
\begin{thm}\label{weakbeta} Let $0<\gamma <1$. Let ${\bf X}=(X_i)_{i\geq 0}$ be a stationary Markov chain with invariant measure $\nu_\gamma$ and transition kernel
$K_\gamma$. There exist two positive constants $C_1(\gamma)$ and
$C_2(\delta, \gamma, k)$ such that, for any $\delta$ in  $]0,
(1-\gamma)/\gamma[$ and any positive integer $k$,
$$
  C_1(\gamma)(n+1)^{\frac{\gamma-1}\gamma}\leq \alpha_k(n) \leq \beta_k(n) \leq C_2(\delta,
  \gamma,k)(n+1)^{\frac{\gamma-1}\gamma+\delta} \, .
$$
\end{thm}

\noindent{\bf Proof of Theorem \ref{weakbeta}.} Applying Proposition
2, Item 2,  in Dedecker and  Prieur (2005), we know that
$$
\nu_\gamma\Big(\sup_{f \in H_{1
,1}}|K_\gamma^nf-\nu_\gamma(f)|\Big)\leq 2\alpha_1(n)\, .
$$
Hence, for any $\varphi$ such that $|\varphi|\leq 1$ and any $f$ in
$H_{1, 1}$,
$$
\nu_\gamma(\varphi \cdot (K_\gamma^nf-\nu_\gamma(f)))=
\nu_\gamma(\varphi \circ T^n \cdot (f-\nu_{\gamma}(f)))\leq
2\alpha_1(n)
$$
The lower bound  for $\alpha_k(n)$ follows from the lower bound for
$\nu_\gamma(\varphi \circ T^n \cdot (f-\nu_{\gamma}(f)))$ given by
Sarig (2002), Corollary 1.

It remains to  prove the upper bound.
 The point is to
approximate the indicator $f_t(x)={\bf 1}_{x \leq t}$ by some
$\alpha$-H\"older function. Let
$$
 f_{t, \epsilon,
 \alpha}(x)=f_t(x)+\Big(1-\Big(\frac{x-t}{\epsilon}\Big)^\alpha\Big)
 {\bf 1}_{t <x \leq t+\epsilon} \, .
$$
This function is $\alpha$-H\"older with H\"older constant
$\epsilon^{-\alpha}$. We now prove the upper bounds  for $k=1$ and
$k=2$ only, the general case being similar. For $k=1$, one has
$$
 K^n(f_{t-\epsilon, \epsilon, \alpha})-\nu_\gamma(f_{t-\epsilon,\epsilon,
 \alpha})-\nu_\gamma([t-\epsilon, t])
   \leq K_\gamma^n(f_t)-\nu_\gamma (f_t) \leq K_\gamma^n(f_{t, \epsilon, \alpha})-\nu_\gamma(f_{t,
  \epsilon, \alpha})+\nu_\gamma([t, t+\epsilon])\, .
$$
Since the  density $g_{\nu_\gamma}$ of $\nu_\gamma$ is such that
$g_{\nu_\gamma}(x) \leq V(\gamma)x^{-\gamma}$, we infer that  for
any real $a$, $\nu_\gamma([a, a+\epsilon])\leq
V(\gamma)\varepsilon^{1-\gamma}(1-\gamma)^{-1}$. Consequently,
$$
|K_\gamma^n(f_t)-\nu_\gamma(f_t)|\leq \epsilon^{-\alpha}\sup_{f\in
H_{\alpha,1}}|K_\gamma^n(f)-\nu_\gamma(f)|  +
\frac{V(\gamma)}{1-\gamma}\epsilon^{1-\gamma} \, .
$$
Applying Theorem \ref{maum} with $k=1$, we obtain that
$$\nu_\gamma\Big(\sup_{t \in
[0,1]}|K_\gamma^n(f_t)-\nu_\gamma(f_t)|\Big) \leq C(\alpha,
1)\epsilon^{-\alpha}(\ln(n+1))^2(n+1)^{\frac{\gamma-1}{\gamma}}+
\frac{V(\gamma)}{1-\gamma}\epsilon^{1-\gamma} \, . $$ The optimal
$\epsilon$ is equal to $$\epsilon= \Big(\frac{\alpha C(\alpha,
1)(\ln(n+1))^2
(n+1)^{\frac{\gamma-1}{\gamma}}}{V(\gamma)}\Big)^{\frac{1}{\alpha+1-\gamma}}\,
.$$ Consequently, for some positive constant $D(\gamma, \alpha)$,
one has
$$
\nu_\gamma\Big(\sup_{t \in
[0,1]}|K_\gamma^n(f_t)-\nu_\gamma(f_t)|\Big)\leq D(\gamma,
\alpha)\Big((\ln(n+1))^2(n+1)^{\frac{\gamma-1}{\gamma}}\Big)^{\frac{1-\gamma}{\alpha+1-\gamma}}\,
.
$$
Choosing $ \alpha<\delta\gamma(1-\gamma)/(1-\gamma(1+\delta)),$ the
result follows for $k=1$.

We now prove the result for $k=2$. Clearly, the four following
inequalities hold:
\begin{eqnarray*}
K_\gamma^n(f_t^{(0)}K_\gamma^m f_s^{(0)})&\leq&
K_\gamma^n(f^{(0)}_{t, \epsilon,
 \alpha}K_\gamma^mf^{(0)}_{s, \epsilon,
 \alpha}) + \nu_\gamma([t, t+ \epsilon]) + \nu_\gamma([s,
 s+\epsilon]) \, , \\
K_\gamma^n(f_t^{(0)}K_\gamma^m f_s^{(0)})&\geq&
K_\gamma^n(f^{(0)}_{t-\epsilon, \epsilon,
 \alpha}K_\gamma^mf^{(0)}_{s-\epsilon, \epsilon,
 \alpha}) - \nu_\gamma([t-\epsilon, t]) - \nu_\gamma([s-\epsilon,
 s]) \, , \\
 \nu_\gamma(f_t^{(0)}K_\gamma^m f_s^{(0)})&\geq& \nu_\gamma(f^{(0)}_{t, \epsilon,
 \alpha}K_\gamma^mf^{(0)}_{s, \epsilon,
 \alpha}) -2 \nu_\gamma([t, t+ \epsilon]) - \nu_\gamma([s,
 s+\epsilon]) \, , \\
\nu_\gamma(f_t^{(0)}K^m f_s^{(0)})&\leq&
\nu_\gamma(f^{(0)}_{t-\epsilon, \epsilon,
 \alpha}K_\gamma^mf^{(0)}_{s-\epsilon, \epsilon,
 \alpha}) +2 \nu_\gamma([t-\epsilon, t]) + \nu_\gamma([s-\epsilon,
 s]) \, .
\end{eqnarray*}
 Consequently,
$$
|K_\gamma^n(f_t^{(0)}K_\gamma^m
f_s^{(0)})-\nu_\gamma(f_t^{(0)}K_\gamma^m f_s^{(0)})| \leq
\epsilon^{-\alpha}\sup_{f,g \in
H_{\alpha,1}}|K_\gamma^n(f^{(0)}K_\gamma^mg^{(0)})-\nu_\gamma(f^{(0)}K_\gamma^mg^{(0)})|
+ \frac{5V(\gamma)}{1-\gamma}\epsilon^{1-\gamma}\, .
$$
Applying Theorem \ref{maum}, we obtain that
$$
\nu_\gamma\Big(\sup_{t \in [0,1]}|K_\gamma^n(f_t^{(0)}K_\gamma^m
f_s^{(0)})-\nu_\gamma(f_t^{(0)}K_\gamma^m f_s^{(0)})|\Big) \leq
C(\alpha,2)\epsilon^{-\alpha}(\ln(n+1))^2(n+1)^{\frac{\gamma-1}{\gamma}}+
\frac{5V(\gamma)}{1-\gamma}\epsilon^{1-\gamma}\, ,
$$
and the proof can be completed as for $k=1$.
\section{Central limit theorems}
 In this section we give a central limit theorem for $S_n(f-\nu_\gamma(f))$ when $f$ belongs to the class
 ${\mathcal C}(M, p, \mu)$  defined in the introduction.
 Note that any function $f$ with bounded variation (BV) such that $|f|\leq M_1$ and $\|df\|\leq M_2$
 belongs to the class ${\mathcal C}(M_1+2M_2, \infty, \mu)$. Hence, any BV function $f$ belongs
to ${\mathcal C}(M, \infty, \mu)$ for some $M$ large enough. If $g$
is monotonic on some open interval of ${\mathbb R}$ and null
elsewhere, and if $\mu(|g|^p) \leq M^p$, then $g$ belongs to
$\mathrm {Mon}(M, p, \mu)$. Conversely, any function in ${\mathcal
C}(M, p, \mu)$ belongs to ${\mathbb L}^q(\mu)$ for $1\leq q<p$.
\begin{thm}\label{CLTgeneral}
Let ${\bf X}=(X_i)_{i\geq 0}$ be a stationary and ergodic (in the
ergodic theoretic sense) Markov chain with invariant measure $\mu$
and transition kernel $K$. Assume that $f$ belongs to ${\mathcal
C}(M, p, \mu)$ for some $M>0$ and some $p\in ]2, \infty]$, and that
$$
  \sum_{k>0}(\alpha_1(k))^{\frac{p-2}{p}} < \infty \, .
$$
The following results hold: \begin{enumerate}
\item The series
$$
\sigma^2(\mu, K, f)=\mu((f-\mu(f))^2)+ 2\sum_{k >0}
\mu((f-\mu(f))K^k(f))
$$
converges to some non negative constant, and
$n^{-1}\mathrm{Var}(\sum_{i=1}^n f(X_i))$ converges to
$\sigma^2(\mu, K, f)$.
\item Let $(D([0, 1], d)$ be the space of cadlag functions from $[0,
1]$ to ${\mathbb R}$ equipped with the Skorohod metric $d$. The
process $\{
 n^{-1/2}\sum_{i=1}^{[nt]} (f(X_i)-\mu(f)), t \in [0, 1]\}$
 converges in distribution in $(D([0, 1], d)$ to $\sigma(\mu, K,
 f)W$, where $W$ is a standard Wiener process.
\item One has the representation  $$f(X_1)-\mu(f)=m(X_1,
X_0)+g(X_1)-g(X_0)$$ with $\mu(|g|^{p/(p-1)})<\infty$, ${\mathbb
E}(m(X_1, X_0)|X_0)=0$ and ${\mathbb E}(m^2(X_1, X_0))=\sigma^2(\mu,
K, f)$.
\end{enumerate}
\end{thm}
\begin{cor}\label{CLTT}
Let $\gamma \in ]0, 1/2[$. If $f$ belongs to the class ${\mathcal
C}(M, p, \nu)$ for some $M>0$ and some $p> (2-2\gamma)/(1-2\gamma)$,
then $n^{-1/2}S_n(f-\nu_\gamma(f))$ converges in distribution to
${\mathcal N}(0, \sigma^2(\nu_\gamma, K_\gamma, f))$.
\end{cor}
\begin{rem}
We infer from Corollary (\ref{CLTT}) that the central limit theorem
holds for any BV function  provided $\gamma<1/2$.  Under the same
condition on $\gamma$, Young (1999) has proved that the central
limit theorem holds for any $\alpha$-H\"older function. For the map
$\theta_\gamma(x)=x(1-x^\gamma)^{-1/\gamma}-[x(1-x^\gamma)^{-1/\gamma}]$
and $\gamma<1/2$, the central limit theorem for BV functions is a
consequence of Corollary 1.7(i) in Raugi (2004).
\end{rem}
\noindent{\bf Two simple examples.} \begin{enumerate} \item Assume
 that $f$ is positive and non increasing on $]0, 1[$, with $f(x)\leq
C x^{-a}$ for some $a\geq 0$. Since the density $g_{\nu_\gamma}$ of
$\nu_{\gamma}$ is such that $g_{\nu_\gamma}(x) \leq
V(\gamma)x^{-\gamma}$, we infer that
$$\nu_{\gamma}(f>t)\leq
\frac{C^{\frac{1-\gamma}{a}}V(\gamma)}{1-\gamma}t^{-\frac{1-\gamma}{a}}\,
.$$ Hence the CLT holds as soon as $a<\frac{1}{2}-\gamma$.
\item
Assume now that $f$ is positive and non decreasing on $]0, 1[$ with
$f(x)\leq C(1-x)^{-a}$ for some $a\geq 0$. Here
$$\nu_\gamma(f>t)\leq
\frac{V(\gamma)}{1-\gamma}\Big(1-\Big(1-\Big(\frac{C}{t}\Big)^{1/a}\Big)^{1-\gamma}\Big)\,
.$$ Hence the CLT holds as soon as
$a<\frac{1}{2}-\frac{\gamma}{2(1-\gamma)}$.
\end{enumerate}

\medskip

 \noindent{\bf Proof of Theorem \ref{CLTgeneral}.}  Let $f$ in ${\mathcal C}(M, p, \mu)$. From Dedecker and Rio
(2000), Items (1) and (2) of Theorem \ref{CLTgeneral} hold as soon
as $$ \sum_{n>0} \|(f(X_0)-\mu(f))({\mathbb
E}(f(X_n)|X_0)-\mu(f))\|_1 < \infty\, .
$$
Assume first that  $f=\sum_{i=1}^k a_i g_i$, where  $\sum_{i=1}^k
|a_i|\leq 1$, and $g_i$ belongs to $\mathrm{Mon}(M, p, \mu)$.
Clearly, the series on left side is bounded by
$$
\sum_{i=1}^k\sum_{j=1}^k|a_ia_j|\sum_{n>0}\|(g_i(X_0)-\mu(g_i))({\mathbb
E}(g_j(X_n)|X_0)-\mu(g_j))\|_1 \, .
$$
Here, we use the following lemma
\begin{lem}\label{deb2}
Let $g_i$ and $g_j$ be two functions in $\mathrm{Mon}(M, p, \mu)$
for some $p\in ]2, \infty]$. For any $1\leq q \leq p$ one has
$$
\|{\mathbb E}(g_j(X_n)|X_0)-\mu(g_j)\|_q \leq
2M\Big(\frac{p}{p-q}\Big)^{1/q}(2\alpha_1(n))^{\frac{p-q}{pq}}\, .
$$
For any $1\leq q<p/2$, one has
$$
\|(g_i(X_0)-\mu(g_i))({\mathbb E}(g_j(X_n)|X_0)-\mu(g_j))\|_q \leq
4M^{2}\Big(\frac{p}{p-2q}\Big)^{1/q}(2\alpha_1(n))^{\frac{p-2q}{pq}}\,
.
$$
\end{lem}
>From Lemma \ref{deb2} with $q=1$, we conclude that
\begin{equation}\label{borne} \sum_{n>0} \|(f(X_0)-\mu(f))({\mathbb
E}(f(X_n)|X_0)-\mu(f))\|_1 \leq
\frac{4pM^{2}}{p-2}\sum_{n>0}(2\alpha_1(n))^{\frac{p-2}{p}}\, .
\end{equation}
Since the bound (\ref{borne}) is true for any function
$f=\sum_{i=1}^k a_i g_i$, it is true also for any $f$ in  ${\mathcal
C}(M, p, \mu)$, and Items (1) and (2) follow.

The last assertion is rather standard.  From the first inequality of
Lemma \ref{deb2} with $q=p/(p-1)$, we  infer that if
$\sum_{n>0}(\alpha_1(n))^{(p-2)/p}<\infty$, then $ \sum_{n>0}
\|{\mathbb E}(f(X_n)|X_0)-\mu(f)\|_{p/(p-1)} <\infty$ for any $f$ in
${\mathcal C}(M, p, \mu)$. It follows that  $g(x)=\sum_{k=1}^\infty
{\mathbb E}(f(X_k)-\mu(f)|X_0=x)$ belongs to ${\mathbb
L}^{p/(p-1)}(\mu)$ and that $m(X_1, X_0)= \sum_{k\geq 1} ({\mathbb
E}(f(X_k)|X_0)-{\mathbb E}(f(X_k)|X_1))$ belongs  to ${\mathbb
L}^{p/(p-1)}$. Clearly
$$
f(X_1)-\mu(f)=m(X_1, X_0)+ g(X_0)-g(X_1) \, ,
$$
with ${\mathbb E}(m(X_1, X_0)|X_0)=0$. Moreover, it follows from the
preceding result that
$$
\lim_{n \rightarrow \infty} \frac {1}{\sqrt n} \Big\|\sum_{k=1}^n
m(X_k, X_{k-1})\Big \|_1=\lim_{n \rightarrow \infty} \frac{1}{\sqrt
n} \Big\|\sum_{k=1}^n (f(X_k)-\mu(f))\Big \|_1 \leq \sigma(\mu, K,
f) \, .
$$
By Theorem 1 in Esseen an Janson (1985), it follows that ${\mathbb
E}(m^2(X_1, X_0))=\sigma^2(\mu, K, f)$.

\medskip

\noindent{\bf Proof of Lemma \ref{deb2}.}  We only prove the second
inequality (the proof of the first one is easier). Let $r=q/(q-1)$
and let $B_r(\sigma(X_0))$ be the set of $\sigma(X_0)$-measurable
random variables such that $\|Y\|_r\leq 1$. By duality,
\begin{eqnarray*}
\|(g_i(X_0)-\mu(g_i))({\mathbb E}(g_j(X_n)|X_0)-\mu(g_j))\|_q&=&
\sup_{Y \in B_r(\sigma(X_0))} {\mathbb
E}(Y(g_i(X_0)-\mu(g_i))(g_j(X_n)-\mu(g_j)))\\
&=& \sup_{Y \in B_r(\sigma(X_0))} \mathrm
{Cov}(Y(g_i(X_0)-\mu(g_i),g_j(X_n))\, .
\end{eqnarray*}
Define the  coefficients $\alpha_{k, g}(n)$ of the sequence
$(g(X_i))_{i \geq 0}$  as in Section \ref{depcoef} with $g\circ f_t$
instead of $f_t$. If $g$ is monotonic on some open interval of
${\mathbb R}$ and null elsewhere, the set $\{x:g(x)\leq t\} $ is
either some interval  or the complement of some interval, so that
$\alpha_{k, g}(n)\leq 2^k\alpha_k (n)$. Let $Q_Y$ be the generalized
inverse of the tail function $t\rightarrow {\mathbb P}(|Y|>t)$. From
Theorem 1.1 and Lemma 2.1 in Rio (2000), one has that
\begin{eqnarray*}
\mathrm {Cov}(Yg_i(X_0),g_j(X_n))&\leq& 2\int_0^{\alpha_{1,
g_i}(n)} Q_Y(u)Q_{g_i(X_0)}(u) Q_{g_j(X_0)}(u) du\\
&\leq& 2\int_0^{2\alpha_{1}(n)} Q_Y(u)Q_{g_i(X_0)}(u)
Q_{g_j(X_0)}(u) du\, .
\end{eqnarray*}
In the same way, applying first Theorem 1.1 in Rio (2000) and next
Fr\'echet's inequality (1957) (see also Inequality (1.11$b$) in Rio
(2000)),
\begin{eqnarray*}
\mathrm {Cov}(Y\mu(g_i),g_j(X_n))&\leq&2
\mu(|g_i|)\int_0^{2\alpha_1(n)} Q_Y(u) Q_{g_j(X_0)}(u) du
\\&\leq& 2\int_0^{2\alpha_1(n)} Q_Y(u)Q_{g_i(X_0)}(u)
Q_{g_j(X_0)}(u) du\, .
\end{eqnarray*}
Since $\int_0^1 Q^r_Y(u) du \leq 1$, it follows that
$$
\|(g_i(X_0)-\mu(g_i))({\mathbb E}(g_j(X_n)|X_0)-\mu(g_j))\|_q \leq 4
\Big(\int_0^{2\alpha_1(n)} Q^q_{g_i(X_0)}(u) Q^q_{g_j(X_0)}(u)
du\Big)^{1/q}\, .
$$
Since $g_i$ and $g_j$ belong to $\mathrm{Mon}(M, p, \mu)$ for some
$p>2q$, we have that $Q_{g_i(X_0)}(u)$ and $Q_{g_j(X_0)}(u)$ are
smaller than $M u^{-1/p}$, and the result follows.

\medskip

\noindent{\bf Proof of Corollary \ref{CLTT}.} We have seen that
$(T_\gamma^1, \ldots, T_\gamma^n)$ is distributed as $(X_n, \ldots,
X_1)$ where $(X_i)_{i \geq 0}$ is the stationary Markov chain with
invariant measure $\nu_\gamma$ and transition kernel $K_\gamma$.
Consequently, on the probability space $([0, 1], \nu_\gamma)$, the
sum $S_n(f-\nu_\gamma(f))$ is distributed as $\sum_{i=1}^n
(f(X_i)-\nu_\gamma(f))$, so that $n^{-1/2}S_n(f-\nu_\gamma(f))$
satisfies the central limit theorem if and only if
$n^{-1/2}\sum_{i=1}^n (f(X_i)-\nu_\gamma(f))$ does. Moreover, we
infer from Theorem \ref{weakbeta} that $$
\alpha_1(n)=O(n^{\frac{\gamma-1}{\gamma}+\epsilon})$$ for any
$\epsilon>0$. Consequently, if $p>(2-2\gamma)/(1-2\gamma)$, one has
that $\sum_{k>0}(\alpha_1(n))^{\frac{p-2}{p}} < \infty$ so that
Theorem \ref{CLTgeneral} applies: the central limit theorem holds
provided that $f$ belongs to ${\mathcal C}(M, p, \nu_\gamma)$.
\section{Rates of convergence in the CLT}
Let $c$ be some concave  function from ${\mathbb R}^+$ to ${\mathbb
R}^+$, with $c(0)=0$. Denote by $\text{Lip}_c$ the set of functions
$g$ such that
$$
|g(x)-g(y)|\leq c(|x-y|) \, .
$$
When $c(x)=x^\alpha$ for $\alpha \in ]0, 1]$, we have $\text
{Lip}_c=H_{\alpha, 1}$. For two probability measures $P, Q$ with
finite first moment, let
$$
d_c(P,Q)=\sup_{g \in \text{Lip}_c}|P(f)-Q(f)| \, .
$$
When $c=$Id, we write $d_c=d_1$. Note that $d_1(P, Q)$ is the
so-called Kantorovi\v c distance between $P$ and $Q$.
\begin{thm}\label{rates}
Let ${\bf X}=(X_i)_{i\geq 0}$ be a stationary  Markov chain with
invariant measure $\mu$ and transition kernel $K$. Let
$\sigma^2(f)=\sigma^2(\mu, K, f)$ be the non-negative number defined
in Theorem \ref{CLTgeneral}, and let $G_{\sigma^2(f)}$ be the
Gaussian distribution with mean 0 and variance $\sigma^2(f)$. Let
$P_n(f)$ be the distribution of the normalized sum
$n^{-1/2}\sum_{i=1}^n (f(X_i) -\mu(f))$.
\begin{enumerate}
\item Assume that $f$ belongs to ${\mathcal C}(M, p, \mu)$ for some $M>0$ and
some $p\in ]2, \infty]$, and that
$$
  \sum_{k>0}(\alpha_{1}(k))^{\frac{p-2}{p}} < \infty \, .
$$
If $\sigma^2(f)=0$, then   $d_c(P_n(f),
\delta_{\{0\}})=O(c(n^{-1/2}))$.
\item If $f$ belongs to ${\mathcal C}(M, p, \mu)$ for some $M>0$ and
some $p\in ]3, \infty]$, and if
$$
  \sum_{k>0}k(\alpha_{3}(k))^{\frac{p-3}{p}} < \infty \, ,
$$
then $d_c(P_n(f), G_{\sigma^2(f)})=O(c(n^{-1/2}))$.
\item If $f$ belongs to ${\mathcal C}(M, p, \mu)$ for some $M>0$ and
some $p\in ]3, \infty]$, and if
$$
  \alpha_{2}(k)=O(k^{-(1+\delta) p/(p-3)}) \quad \text{for some
  $\delta \in ]0, 1[$},
$$
then $d_c(P_n(f), G_{\sigma^2(f)})=O(c(n^{-\delta/2}))$.
\end{enumerate}
\end{thm}
\begin{cor}\label{Krov}
Let $\delta \in ]0, 1]$ and  $\gamma <1/(2+\delta)$, and let
$\mu_n(f)$ be the distribution of $n^{-1/2}S_n(f-\nu_\gamma(f))$. If
$f$ belongs to the class ${\mathcal C}(M, p, \nu_\gamma)$ for some
$M>0$ and some $p>(3-3\gamma)/(1-(2+\delta)\gamma)$, then
$d_c(\mu_n(f), G_{\sigma^2(f)})=O(c(n^{-\delta/2}))$, where
$\sigma^2(f)=\sigma^2(\nu_\gamma, K_\gamma, f)$.
\end{cor}
\begin{rem} We infer from Corollary \ref{Krov} that if $f$ is $BV$, then $d_1(\mu_n(f), G_{\sigma^2(f)})=O(n^{-1/2})$ if
$\gamma<1/3$, and $d_1(\mu_n(f),G_{\sigma^2(f)})=O(n^{-\delta/2})$
if $\gamma<1/(2+\delta)$. Denote by $d_{BV}(P, Q)$ the uniform
distance between the distribution functions of $P$ and $Q$. If $f$
is $\alpha$-H\"older, Gou\"ezel (2005, Theorem 1.5) has proved that
$d_{BV}(\mu_n(f), G_{\sigma^2}(f))=O(n^{-1/2})$ if $\gamma<1/3$, and
$d_{BV}(\mu_n(f),G_{\sigma^2(f)})=O(n^{-\delta/2})$ if
$\gamma=1/(2+\delta)$. In fact, from a general result of Bolthausen
(1982) for Harris recurrent Markov chains, we conjecture that the
results of Corollary \ref{Krov} are true with $d_{BV}$ instead of
$d_1$.
\end{rem}
\noindent{\bf Two simple examples (continued).} \begin{enumerate}
\item
Assume
 that $f$ is positive and non increasing on $[0, 1]$, with $f(x)\leq
C x^{-a}$ for some $a\geq 0$. Let $\delta \in ]0, 1]$ and $\gamma
<1/(2 + \delta)$. If $a<\frac{1}{3}-\frac{(2+\delta)\gamma}{3}$,
then $d_c(\mu_n(f), G_{\sigma^2(f)})=O(c(n^{-\delta/2}))$.
\item Assume
 that $f$ is positive and non increasing on $[0, 1]$, with $f(x)\leq
C (1-x)^{-a}$ for some $a\geq 0$. Let $\delta \in ]0, 1]$ and
$\gamma <1/(2 + \delta)$. If
$a<\frac{1}{3}-\frac{(1+\delta)\gamma}{3(1-\gamma)}$, then
$d_c(\mu_n(f), G_{\sigma^2(f)})=O(c(n^{-\delta/2}))$.
\end{enumerate}

\medskip

\noindent{\bf Proof of Theorem \ref{rates}.} From the Kantorovi\v
c-Rubin\v ste\u\i n  theorem (1957), there exists a probability
measure $\pi$ with margins $P$ and $Q$, such that $ d_1(P, Q)=\int
|x-y| \pi(dx, dy) $. Since $c$ is concave, we then have
$$
d_c(P, Q)=\sup_{f \in H_c}\Big|\int (f(x)-f(y))\pi(dx, dy)\Big| \leq
\int c(|x-y|) \pi(dx, dy) \leq c(d_1(P, Q)) \, .
$$
Hence, it is enough to prove the theorem for $d_1$ only.

If  $\sum_{k>0}(\alpha_{1}(k))^{(p-2)/p} < \infty$, $f$ belongs to
${\mathcal C}(M, p, \mu)$ for some $M>0$ and some $p\in ]2,
\infty]$,  and $\sigma^2(f)=0$, it follows from Theorem
\ref{CLTgeneral} that $f(X_1)=g(X_0)-g(X_1)$ with $\mu(|g|)<
\infty$. Hence
$$
d_1(P_n(f), \delta_{\{0\}}) \leq \frac{2\mu(|g|)}{\sqrt n} \, ,
$$
and Item (1) is proved.

>From now, we assume that $\sigma^2(f)>0$ (otherwise, the result
follows from Item (1)). If  $f=g_1-g_2$,   where $g_1, g_2$ belong
to $\mathrm {Mon}(M, p, \mu)$ for some
 $M>0$ and
some $p\in ]3, \infty]$,
 Item (2) of Theorem \ref{rates} follows from Theorem 3.1(b) in Dedecker and Rio
 (2007). In fact the proof remains unchanged if $f$ belongs to
 ${\mathcal C}(M, p, \mu)$ for some $M>0$ and
some $p\in ]3, \infty]$.

It remains to prove Item (3).  Let $Y_k=f(X_k)-\mu(f)$,
$\sigma^2(f)=\sigma^2$,  and $s_m=\sum_{i=1}^m Y_i$. Define
$$
W_m=A_m+B_m,\quad \text{with} \quad  A_m={\mathbb E}( s_m^2| X_0)-m
\sigma^2 \quad \text{and} \quad B_m= 2\sum_{k=1}^m{\mathbb
E}\Big(Y_k \sum_{i>m} Y_i\Big|X_0\Big)\, .
$$
>From  Theorem 2.2 in Dedecker and Rio (2007), we have that, if
$\sum_{k>0}\|Y_0{\mathbb E}(Y_k|X_0)\|_1 < \infty$,
\begin{equation}\label{deri}
\sqrt n d_1(P_n(f),G_{\sigma^2}) \leq C\ln(n)+
\sum_{m=1}^{[\sqrt{2n}]}\frac{\|(|Y_0|+2\sigma)W_m\|_1}{m\sigma^2}
+D_{1,n}+ D_{2,n}\, ,
\end{equation}
where
$$
D_{1,n}=\sum_{m=1}^n {\frac{1}{\sigma\sqrt{m}}}   \sum_{i\geq m}
\|Y_0{\mathbb E}(Y_i|X_0)\|_1  \quad \text{and} \quad
D_{2,n}=\sum_{m=1}^n {\frac{1}{2\sigma^2m}}
 \sum_{k=1}^m \|(\sigma^2+ Y_0^2) {\mathbb E}(Y_k|X_0)\|_1  .
$$
>From Lemma \ref{deb2} with $q=1$, the bound (\ref{borne}) holds for
any $f$ in ${\mathcal C}(M, p, \mu)$ for $p>2$. Consequently, if
$\alpha_{2}(k)=O(k^{-(1+\delta) p/(p-3)})$ for some $\delta \in ]0,
1[$ and $p>3$, then $\sum_{k>0}\|Y_0{\mathbb E}(Y_k|X_0)\|_1 <
\infty$, so that the bound (\ref{deri}) holds. Moreover
$n^{-1/2}D_{1, n} = O(n^{-1/2}\ln(n)\vee n^{-\delta})$. Arguing as
in Lemma \ref{deb2}, one can prove that
$$
\|Y_0^2{\mathbb E}(Y_k|X_0)\|_1 \leq C(M, p)
(\alpha_1(k))^{\frac{p-3}{p}}\, ,
$$
so that $n^{-1/2} D_{2, n}=O(n^{-1/2} \ln(n))$.

Arguing as in Lemma \ref{deb2}, one can prove that, for $0<k<i$,
\begin{equation}\label{petiteborne}
 \|(|Y_0|+2\sigma){\mathbb E}(Y_k Y_i|X_0)\|_1\leq  \|(|Y_0|+2\sigma) Y_k{\mathbb E}(Y_i|X_k)\|_1\leq C(M, p, \sigma)
  (\alpha_1(i-k))^{\frac{p-3}{p}}\, .
\end{equation}
Consequently,
$$
\frac {1}{ \sqrt n}
\sum_{m=1}^{[\sqrt{2n}]}\frac{\|(|Y_0|+2\sigma)B_m\|_1}{m\sigma^2}=
O\Big( \frac {1}{ \sqrt
n}\sum_{m=1}^{[\sqrt{2n}]}\frac{1}{m\sigma^2}\sum_{k=1}^m
\sum_{i>m}\frac{1}{(i-k)^{1+\delta}}\Big)=O(n^{-\delta/2})\, .
$$

Now,
$$
\frac{\|(|Y_0|+2\sigma)A_m\|_1}{m} \leq \frac{2}{m}\sum_{i=1}^m
\sum_{j=i}^m \|(|Y_0|+2\sigma)({\mathbb E}(Y_iY_j|X_0)-{\mathbb
E}(Y_iY_j))\|_1 +(\|Y_0\|_1+2\sigma)\Big|\frac1m {\mathbb
E}(s_m^2)-\sigma^2\Big|\, .
$$
For the second term on right hand, we have
$$
\Big|\frac1m {\mathbb E}(s_m^2)-\sigma^2\Big|\leq
2\sum_{k=1}^\infty\frac {k\wedge m}m |{\mathbb E}(Y_0 Y_k)|=O\Big(
\sum_{k>0}\frac {k\wedge m}m
(\alpha_1(k))^{\frac{p-2}{p}}\Big)=O(m^{-\delta})\, ,
$$
so that
$$
\frac{1}{\sqrt n}\sum_{m=1}^{[\sqrt{2n}]} \Big|\frac1m {\mathbb
E}(s_m^2)-\sigma^2\Big|=O(n^{-\delta/2}) \, .
$$

To complete the proof of the theorem, it remains to prove that
\begin{equation}\label{ouf}
\frac{1}{\sqrt n}\sum_{m=1}^{[\sqrt{2n}]} \frac{2}{m}\sum_{i=1}^m
\sum_{j=i}^m \|(|Y_0|+2\sigma)({\mathbb E}(Y_iY_j|X_0)-{\mathbb
E}(Y_iY_j))\|_1 =O(n^{-\delta/2})\, .
\end{equation}
Applying first (\ref{petiteborne}), we have for $j>i$,
\begin{equation}\label{premborn}
\|(|Y_0|+2\sigma)({\mathbb E}(Y_iY_j|X_0)-{\mathbb
E}(Y_iY_j))\|_1\leq 2C(M, p, \sigma)
  (\alpha_1(j-i))^{\frac{p-3}{p}}\, .
\end{equation}
We need a second bound for this quantity. Assume first that
$f=\sum_{i=1}^k a_i g_i$, where $\sum_{i=1}^k|a_i|\leq 1$ and $g_i$
belongs to $\mathrm{Mon}(M,p,\mu)$. Let $g_i^{(0)}=g_i-\mu(g_i)$. We
have that
\begin{multline*}
\|Y_0({\mathbb E}(Y_iY_j|X_0)-{\mathbb E}(Y_iY_j))\|_1\\
\leq \sum_{l=1}^k\sum_{q=1}^k \sum_{r=1}^k
|a_la_qa_r|\|g^{(0)}_l(X_0)({\mathbb
E}(g^{(0)}_q(X_i)g^{(0)}_r(X_j)|X_0)-{\mathbb
E}(g^{(0)}_q(X_i)g^{(0)}_r(X_j)))\|_1\, .
\end{multline*}
For three real-valued random variables $A, B, C$, define the numbers
$\bar{\alpha}(A, B)$ and  $\bar{\alpha}(A, B, C)$ by
\begin{eqnarray*}
\bar{\alpha}(A, B)&=& \sup_{s, t\in {\mathbb R}}|\mathrm{Cov}({\bf
1}_{A\leq s},{\bf 1}_{B\leq t})|\\ \bar{\alpha}(A, B, C)&=&\sup_{s,
t, u \in {\mathbb R}}|{\mathbb E}(({\bf 1}_{A\leq s}-{\mathbb
P}(A\leq s))({\bf 1}_{B\leq t}-{\mathbb P}(B\leq t))({\bf 1}_{C\leq
u}-{\mathbb P}(C\leq u)))|\,
\end{eqnarray*}
(note that $\bar{\alpha}(A, B, B)\leq \bar{\alpha}(A, B)$). Let
$$A=|g^{(0)}_l(X_0)|\mathrm{sign}\{{\mathbb
E}(g^{(0)}_q(X_i)g^{(0)}_r(X_j)|X_0)-{\mathbb
E}(g^{(0)}_q(X_i)g^{(0)}_r(X_j))\}\, , $$ and note that $Q_{A}=
Q_{g^{(0)}_l(X_0)}$. From Proposition 6.1 and Lemma 6.1 in Dedecker
and Rio (2007), we have that
\begin{multline*}
\|g^{(0)}_l(X_0)({\mathbb
E}(g^{(0)}_q(X_i)g^{(0)}_r(X_j)|X_0)-{\mathbb
E}(g^{(0)}_q(X_i)g^{(0)}_r(X_j)))\|_1={\mathbb E}((A-{\mathbb E}(A))g^{(0)}_q(X_i)g^{(0)}_r(X_j))\\
\leq 16 \int_0^{\bar{\alpha}(A, g_q(X_i),
g_r(X_j))/2}Q_{g^{(0)}_l(X_0)}(u)Q_{g_q(X_0)}(u)Q_{g_r(X_0)}(u)du\,
.
\end{multline*}
Note that $Q_{g^{(0)}_l(X_0)}\leq Q_{g_l(X_0)}+\|g_l(X_0)\|_1$.
Hence, by Fr\'echet's inequality (1957),
\begin{multline*}
 \int_0^{\bar{\alpha}(A, g_q(X_i),
g_r(X_j))/2}Q_{g^{(0)}_l(X_0)}(u)Q_{g_q(X_0)}(u)Q_{g_r(X_0)}(u)du \\
\leq 2 \int_0^{\bar{\alpha}(A, g_q(X_i),
g_r(X_j))/2}Q_{g_l(X_0)}(u)Q_{g_q(X_0)}(u)Q_{g_r(X_0)}(u)du\, .
\end{multline*}
Since $\{g_i(x)\leq t\}$ is some interval of ${\mathbb R}$, we have
that for $j>i\geq 1$
$$\bar{\alpha}(A, g_q(X_i), g_r(X_j))\leq 4 \bar{\alpha}(A, X_i,
X_j)\leq 4\alpha_2(i)\, ,$$ and for $i=j$,
$$\bar{\alpha}(A, g_q(X_i), g_r(X_i))\leq 4 \bar{\alpha}(A,
X_i, X_i)\leq 4 \bar{\alpha}(X_0, X_i)\leq 4\alpha_1(i)\leq
4\alpha_2(i)\, .$$ Since $Q_{g_i(X_0)}(u)\leq Mu^{-1/p}$, it follows
that, for $1\leq i\leq j$,
$$
\|g_l(X_0)({\mathbb E}(g_q(X_i)g_r(X_j)|X_0)-{\mathbb
E}(g_q(X_i)g_r(X_j)))\|_1\leq \frac{32
M^3p}{p-3}(2\alpha_2(i))^{\frac{p-3}{p}}\, .
$$
Consequently, for any $f$ in ${\mathcal C}(M, p, \mu)$ with $p>3$,
$$
\|Y_0({\mathbb E}(Y_iY_j|X_0)-{\mathbb E}(Y_iY_j))\|_1\leq \frac{32
M^3p}{p-3}(2\alpha_2(i))^{\frac{p-3}{p}}\, .
$$
In the same way,
$$
2\sigma\|{\mathbb E}(Y_iY_j|X_0)-{\mathbb E}(Y_iY_j)\|_1\leq
\frac{32\sigma M^2p}{p-2}(2\alpha_2(i))^{\frac{p-2}{p}}\, .
$$
It follows that, for any $1\leq i\leq j$,
\begin{equation}\label{deuxborne}
\|(|Y_0|+2\sigma)({\mathbb E}(Y_iY_j|X_0)-{\mathbb
E}(Y_iY_j))\|_1\leq D(M, p, \sigma)
  (\alpha_2(i))^{\frac{p-3}{p}}\, .
\end{equation}
Combining (\ref{premborn}) and (\ref{deuxborne}), we infer that
$$
\sum_{i=1}^m \sum_{j=i}^m \|(|Y_0|+2\sigma)({\mathbb
E}(Y_iY_j|X_0)-{\mathbb E}(Y_iY_j))\|_1 =O(m^{1-\delta})\, ,
$$
and (\ref{ouf}) easily follows. This completes the proof.
\section{Moment inequalities}
\begin{thm}\label{M}
Let ${\bf X}=(X_i)_{i\geq 0}$ be a stationary  Markov chain with
invariant measure $\mu$ and transition kernel $K$. If $f$ belong to
${\mathcal C}(M, p, \mu)$ for some $M>0$ and some $p>2$, then, for
any $2\leq q <p$
$$
\Big \| \sum_{i=1}^n (f(X_i)-\mu(f))\Big\|_q \leq \sqrt {2q}\Big(n
\|f(X_0)-\mu(f)\|_q^{2}+4M^{2}\Big(\frac{p}{p-q}\Big)^{\frac2q}\,
\sum_{k=1}^{n-1}(n-k)(2\alpha_1(k))^{\frac{2(p-q)}{pq}}\Big)^{\frac12}\,.
$$
\end{thm}
\begin{cor}\label{melnic}
Let $0<\gamma<1$. Let $f$ belong to ${\mathcal C}(M, p, \mu)$ for
some $M>0$ and some $p>2$, and let $2\leq q <p$.
\begin{enumerate}
\item
If $\gamma< 2(p-q)/(2(p-q)+pq)$, then
$
  \|S_n(f-\nu_\gamma(f))\|_q =O(\sqrt n)
$\, .
\item If $2(p-q)/(2(p-q)+pq)\leq\gamma<1$, then, for any
$\epsilon>0$,
$$
  \|S_n(f-\nu_\gamma(f))\|_q =O\Big( n^{1+\epsilon-\frac{(1-\gamma)(p-q)}{\gamma p
  q}}\Big)\, .
$$
\end{enumerate}
\end{cor}
\begin{rem}
Assume that $\gamma<(p-2)/(2p-2)$. By Chebichev inequality applied
with $2\leq q< 2p(1-\gamma)/(\gamma p +2(1 -\gamma))$, we infer from
Item (1) that for any $\epsilon>0$,
\begin{equation*}\label{MN}
  \nu_\gamma\Big(\frac{1}{n}|S_n(f-\nu_\gamma(f))|>x\Big)\leq
  \frac{C}{(nx^2)^{p(1-\gamma)/(\gamma p+2(1-\gamma))-\epsilon}}\, .
\end{equation*}
Assume now that $(p-2)/(2p-2)\leq \gamma < 1$. By Chebichev
inequality applied with $q=2$, we infer from Item (2) that for any
$\epsilon>0$,
\begin{equation*}\label{MN}
  \nu_\gamma\Big(\frac{1}{n}|S_n(f-\nu_\gamma(f))|>x\Big)\leq
  \frac{C}{x^2n^{(p-2)(1-\gamma)/\gamma p-\epsilon}}\, .
\end{equation*}
When $f$ is BV (case $p=\infty$) and $\gamma<1$, we obtain that, for
any $\epsilon>0$ and any $x>0$,
\begin{equation*}\label{MN}
  \nu_\gamma\Big(\frac{1}{n}|S_n(f-\nu_\gamma(f))|>x\Big)\leq
  \frac{C(x)}{n^{(1-\gamma)/\gamma -\epsilon}}\, .
\end{equation*}
Note that Melbourne and Nicol (2007) obtained the same bound when
$f$ is $\alpha$-H\"older and $\gamma<1/2$.
\end{rem}

\noindent{\bf Two simple examples (continued).} \begin{enumerate}
\item
Assume
 that $f$ is positive and non increasing on $[0, 1]$, with $f(x)\leq
C x^{-a}$ for some $a>0$. If  $a<\frac 12 -\gamma$ and $2\leq q
<\frac{2(1-\gamma)}{\gamma+2a}$, then $\|S_n(f-\nu_\gamma(f))\|_q=
O(\sqrt n)$. If now $a<\frac {1 -\gamma}{2}$ and $2\vee
\frac{2(1-\gamma)}{\gamma+2a}\leq q < \frac{1-\gamma}{a}$, then, for
any $\epsilon>0$,  $$\|S_n(f-\nu_\gamma(f))\|_q=
O\Big(n^{1+\epsilon-\frac{(1-\gamma - aq)}{\gamma
  q}}\Big)\,.
  $$
\item Assume
 that $f$ is positive and non increasing on $[0, 1]$, with $f(x)\leq
C (1-x)^{-a}$ for some $a\geq 0$. If $a<
\frac{1-2\gamma}{2(1-\gamma)}$ and $2\leq q
<\frac{2(1-\gamma)}{\gamma+(1-\gamma)2a}$, then
$\|S_n(f-\nu_\gamma(f))\|_q= O(\sqrt n)$. If $a<\frac 12$ and $2\vee
\frac{2(1-\gamma)}{\gamma+(1-\gamma)2a}\leq q < \frac{1}{a}$, then,
for any $\epsilon>0$,  $$\|S_n(f-\nu_\gamma(f))\|_q=
O\Big(n^{1+\epsilon-\frac{(1-\gamma)(1 - aq)}{\gamma
  q}}\Big)\,.
  $$
\end{enumerate}

\noindent{\bf Proof of Theorem \ref{M}.} From Proposition 4  in
Dedecker and Doukhan (2003) (see also Theorem 2.5 in Rio (2000)), we
have that, for any $q\geq 2$,
$$
\Big \| \sum_{i=1}^n (f(X_i)-\mu(f))\Big\|_q \leq \sqrt{2q} \Big(n
\|f(X_0)-\mu(f)\|_q^{2}+\sum_{k=1}^{n-1}
(n-k)\|(f(X_0)-\mu(f))({\mathbb
E}(f(X_k)|X_0)-\mu(f))\|_{\frac{q}{2}} \Big)^{\frac 12} .
$$
Assume first that  $f=\sum_{i=1}^k a_i g_i$, where  $\sum_{i=1}^k
|a_i|\leq 1$, and $g_i$ belongs to $\mathrm{Mon}(M, p, \mu)$.
Clearly
$$
\|(f(X_0)-\mu(f))({\mathbb E}(f(X_n)|X_0)-\mu(f))\|_{q/2}\leq
\sum_{i=1}^k\sum_{j=1}^k|a_ia_j|\|(g_i(X_0)-\mu(g_i))({\mathbb
E}(g_j(X_n)|X_0)-\mu(g_j))\|_{q/2} \, .
$$
Applying Lemma \ref{deb2}, we obtain that
$$
\|(f(X_0)-\mu(f))({\mathbb E}(f(X_n)|X_0)-\mu(f))\|_{q/2}\leq
4M^{2}\Big(\frac{p}{p-q}\Big)^{2/q}(2\alpha_1(n))^{\frac{2(p-q)}{pq}}\,
.
$$
Clearly, this inequality remains valid for any $f$ in ${\mathcal
C}(M, p, \mu)$, and the result follows.
\section{The empirical distribution function}
\begin{thm}\label{emp}
Let ${\bf X}=(X_i)_{i\geq 0}$ be a stationary  Markov chain with
invariant measure $\mu$ and transition kernel $K$. Let
$F_n(t)=n^{-1} \sum_{i=1}^n {\bf 1}_{X_i \leq t}$  and
$F_\mu(t)=\mu(]-\infty, t])$.
\begin{enumerate}
\item If ${\bf X}$ is ergodic (in the ergodic theoretic sense) and
if $\sum_{k>0} \beta_1(k)<\infty$, then, for any probability $\pi$
on ${\mathbb R}$, the process $\{\sqrt n (F_n(t)-F_\mu(t)), t \in
{\mathbb R}\}$ converges in distribution in ${\mathbb L}^2(\pi)$ to
a tight Gaussian process $G$ with covariance function
\begin{equation*}\label{covemp}
\mathrm{Cov}(G(s), G(t))=C_{\mu, K}(s,t)=\mu(f^{(0)}_t
f^{(0)}_s)+2\sum_{k>0}\mu(f^{(0)}_t K^k f^{(0)}_s)\, .
\end{equation*}
\item Let $(D({\mathbb R}), d)$ be the space of cadlag functions equipped with the Skorohod metric $d$.
 If $\beta_2(k)=O(k^{-2 -\epsilon})$ for some $\epsilon>0$, then the process $\{\sqrt n (F_n(t)-F_\mu(t)), t \in
{\mathbb R}\}$ converges in distribution in $(D({\mathbb R}), d)$ to
a tight Gaussian process $G$ with covariance function $C_{\mu, K}$.
\end{enumerate}
\end{thm}
 \begin{cor}\label{empT} Let $F_{n, \gamma}(t)= n^{-1} \sum_{i=1}^n {\bf 1}_{T_\gamma^i \leq
 t}$.
\begin{enumerate}
\item If $0<\gamma < 1/2$, then,  for any probability $\pi$
on $[0, 1]$, the process $\{ \sqrt{n} (F_{n,
\gamma}(t)-F_{\nu_\gamma}(t)), t \in {[0, 1]}\}$ converges in
distribution in ${\mathbb L}^2(\pi)$ to a tight Gaussian process
$G_\gamma$ with covariance function $C_{\nu_\gamma, K_\gamma}$.
\item If $0<\gamma < 1/3$, the process $\{ \sqrt{n} (F_{n,
\gamma}(t)-F_{\nu_\gamma}(t)), t \in {[0, 1]}\}$ converges in
distribution in $(D([0, 1]), d)$ to a tight Gaussian process
$G_\gamma$ with covariance function $C_{\nu_\gamma, K_\gamma}$.
\end{enumerate}
\end{cor}

\begin{rem}
Denote by $\|\cdot\|_{p, \pi}$ the ${\mathbb L}^p(\pi)$-norm. If
$\gamma< 1/2$, we have that, for any $1\leq p \leq 2$,
\begin{equation}\label{norm}
  {\sqrt n}\|F_{n, \gamma}-F_{\nu_\gamma}\|_{p, \pi} \quad \text{converges
  in distribution to} \quad \|G_\gamma\|_{p, \pi}\, .
\end{equation}
In particular, if $\pi=\lambda$ is the Lebesgue measure on $[0, 1]$
and $q=p/(p-1)$,  we obtain that
$$\frac{1}{\sqrt n}\sup_{\|f'\|_q\leq 1}|S_n(f-\nu_\gamma(f))| \quad
\text{converges in distribution to} \quad \|G_\gamma\|_{p,
\lambda}
\, .$$ For $p=1$
and $q=\infty$, we obtain the limit distribution of  the Kantorovi\v
c distance $d_1(F_{n, \gamma}, F_{\nu_\gamma})$:
$$\sqrt n d_1(F_{n, \gamma},
F_{\nu_\gamma})=\frac{1}{\sqrt n}\sup_{f \in H_{1,
1}}|S_n(f-\nu_\gamma(f))| \quad \text{converges in distribution to}
\quad  \int_0^1 |G_\gamma(t)| dt\, .$$  Now if $\gamma<1/3$, the
limit in (\ref{norm}) holds for any $p\geq 1$.

Note that, for Harris recurrent Markov chains, Item (2) of Theorem
\ref{emp} holds as soon as the sum of the $\beta$-mixing
coefficients of the chain is finite. Hence, we conjecture that Item
(2) of Corollary \ref{empT} remains true for $\gamma<1/2$.

\end{rem}

\noindent{\bf Proof of Theorem \ref{emp}.} Item (1) has been proved
in Dedecker and Merlev\`ede (2007, Theorem 2, Item 2) and Item (2)
in Dedecker and Prieur (2007, Proposition 2).

\medskip

\noindent{\bf Acknowledgments.} Many thanks to Jean-Ren\'e
Chazottes, who pointed out the references to Conze and Raugi (2003)
and Raugi (2004).

\end{document}